# A formula for the n'th digit of $\pi$ and $\pi^n$

Simon Plouffe
January 29, 2022

Abstract

By using an asymptotic formula known for the numbers of Euler and Bernoulli it is possible to obtain an explicit expression of the nth digit of π in decimal or in binary, it also makes it possible to obtain the nth digit of $\pi^n$. The calculation is made from the two inequalities

$$\frac{2(2n)!}{(2\pi)^{2n}}\left(\frac{1}{1-2^{1-2n}}\right) > (-1)^{n+1}B_{2n} > \frac{2(2n)!}{(2\pi)^{2n}}, n = 1, 2, \ldots$$

et

$$\frac{4^{n+1}(2n)!}{\pi^{2n+1}} > E_{2n} > \frac{4^{n+1}(2n)!}{\pi^{2n+1}}\left(\frac{1}{1+3^{-1-2n}}\right), n = 0, 1, \ldots$$

By isolating π in both cases, we can derive an approximation of the latter. The approximation is so good that in fact it allows to extract the nth bit of π, the nth decimal of π and even the nth decimal of de $\pi^n$, n ∈ $\mathbb{Z}^*$.

## Bernoulli numbers and $\pi$

The 2 inequalities are taken from [3], p. 809. By isolating therefore, the number π in the first one obtains:

$$\pi \approx \left(\frac{2n!}{B_n 2^n}\right)^{1/n}$$

Here the $B_n$ are positive and n is even. It is easy to verify that with n = 1000, the error on π is of the order of $0.293193 \times 10^{-303}$, which is less than $2^{-1000}$. So the thousandth position of this expression is the 1000th bit of the number π. As soon as n = 10, the formula is valid: it allows to have the nth position. Note that given the parity of Bernoulli numbers, we obtain 2 bits at each n, that of rand n and that of rank n-1. This is valid from 10 because Bernoulli numbers are a little chaotic when n is small, the rapid growth of these numbers is felt from n = 12.

The good question that comes next is: is it possible to calculate the nth decimal place? The answer is actually yes. Here's how.

Given the known relationship with the Zeta function, it is possible to stretch the formula by adding enough terms to it to have an error of less than $10^{-n}$. In the book of Abramowitz and Stegun [3] precisely in the same chapter we have the formula

$$\zeta(2n) = \frac{(2\pi)^{2n}}{2(2n)!}|B_{2n}|$$

As we know, ζ(2n) can be represented as an infinite product due to Euler, by carefully taking the first 4 terms we therefore add them to our expression to have:

$$\pi \approx \left(\frac{2n!}{B_n 2^n \left(1-\frac{1}{2^n}\right)\left(1-\frac{1}{3^n}\right)\left(1-\frac{1}{5^n}\right)\left(1-\frac{1}{7^n}\right)}\right)^{1/n}$$

Knowing that the next term would be $\left(1-\frac{1}{11^n}\right)$ this assures us that the error committed will be less than 10^(-n). Indeed, by checking with n = 1000 (and even 10000) we find an error of the order of $0.1271934403 \times 10^{-1043}$. We conclude that for each n in the expression it seems possible to obtain the nth decimal of the number π. The difficulty is then to be able to calculate B_2n when n is large. The current record for Bernoulli numbers is 100 million. So, we can calculate the 100,000,000th decimal of π with this process. Note that the calculation of Bernoulli numbers can be done in several ways, one of which requires knowing π with good precision. Until 2002, the Bernoulli calculation record used precisely this formula and some values of the series for ζ(n), see [5] on this subject. Currently, the calculation uses properties of congruences and the Chinese remainder theorem.

The formula is precise enough to get even more if we raise to the power n. In fact we will have

$$\pi^n \approx \left(\frac{2n!}{B_n 2^n \left(1-\frac{1}{2^n}\right)\left(1-\frac{1}{3^n}\right)\left(1-\frac{1}{5^n}\right)\left(1-\frac{1}{7^n}\right)}\right)$$

By checking with n = 1000 we obtain the first 1000 valid digits of π^1000, so the 1000th position if we want. To get the 1000th position of x, just calculate

$$position(n, x) = [10\,\{10^n x\}]$$

[ ] and { } being the integer and fractional parts. If calculating π^n, the formula must be converted to treat the number as a string.

## Euler numbers and $\pi$

We take the expression of the book of A & S page 809, we isolate π and we obtain:

$$\pi \approx \left( \frac{(2n)!\, 2^{2n+2}}{E_{2n}} \right)^{1/(2n+1)}$$

Here the $E_{2n}$ are the positive Euler numbers obtained by the series expansion of $1/cos(x)$. It is easy to verify that with n = 1000 we obtain a precision of 957 decimal digits. But knowing the error made and the first formula, we can deduce that this error is proportional to 9^(-1000). So, we can say that in base 8 the expression gives the nth digit of π in octal, a fortiori in binary.

To go further as with the Bernoulli numbers, it is enough to add only 1 term to have an error lower than 10^(-n). Which gives us

$$\pi \approx \left( \frac{(2n)!\, 2^{2n+2}}{E_{2n}} \left(1 - \frac{1}{3^{2n+1}}\right) \right)^{1/(2n+1)}$$

And if we check with n >>1, for n = 1000 we have here an error on π of the order of 10^(-1198), so the nth decimal of π can be calculated with the formula using the position function (n,x).

To get the nth digit in base 10 of $\pi^n$ using Euler numbers, we add terms until the error is smaller than $10^{-n}$. We have to take more terms here since the powers of π take us away from the decimal point. So we finally have the following formula

$$\pi^{2n+1} \approx \frac{(2n)!\, 2^{2n+2}}{E_{2n}} \left(1 - \frac{1}{3^{2n+1}}\right)\left(1 + \frac{1}{5^{2n+1}}\right)\left(1 - \frac{1}{7^{2n+1}}\right)\left(1 - \frac{1}{9^{2n+1}}\right)$$

For 1/ π, in base 10 still with Euler, we can simply invert

$$\frac{1}{\pi} \approx \left( \frac{E_{2n}}{(2n)!\, 2^{2n+2}} \left(1 - \frac{1}{3^{2n+1}}\right) \right)^{1/(2n+1)}$$

It therefore makes it possible to obtain the nth decimal of 1/ π in base 10. The scope of this formula is limited by the ability to evaluate E_2n . With an Intel icore 9900K clocked at 5 ghz, I managed without too much difficulty to calculate E_(2,300,000). So position 2,300,000 can be calculated by this process.

In the end, we can therefore calculate the nth decimal digit of π and π^n. In the end, we can therefore calculate the nth decimal digit of π and π^n. The same process also makes it possible to calculate the nth decimal of the following numbers.

$$e^{\pi},\ \ln(\pi),\ \sqrt{\pi}$$

By isolating n! from the equation we obtain a much better formula than that of Stirling.

$$n! \approx \frac{(2\pi)^n B_n}{2}$$

The same with Euler numbers

$$n! \approx \frac{\pi^{n+1} E_n}{2^{n+2}}$$

The 2 expressions are in absolute value and n even. A summary analysis of the error shows for $n = 1000$ shows that for the first the error is of the order of $2^{-1000}$, i.e. $10^{-301}$ and $3^{-1000}$ or $10^{-477}$ for the 2nd. With Stirling's formula, the error is of the order of $10^{-22}$ only. As with the previous formulas, we can improve even more if we add terms drawn from the function ζ(n).

But there is better, by comparing two successive terms of the expression with the Bernoulli numbers it allows to get rid of the exponent n. If we make the ratio, in fact we obtain for $\pi^2$ a closed expression.

We had

$$\pi \approx \left(\frac{2n!}{B_n 2^n}\right)^{1/n}$$

By comparing 2 successive terms we get

$$\pi^2 \approx \frac{1}{2} \frac{B_{2n}\,(n+1)(2n+1)}{B_{2n+2}}$$

Which is valid for base 4 given that the error is of the order of $4^{-n}$, by applying the same technique with Euler numbers we then have:

$$\pi^2 \approx \frac{8 E_{2n}\,(n+1)(2n+1)}{E_{2n+2}}$$

Whose error with respect to π^2 is $9^{-n}$. So in base 9. We could use the same arguments for base 10 in both cases by adding terms. Hence the age-old question: is there an interpretation of successive Bernoulli numbers allowing a combinatorial interpretation of the binary digits of $\pi^2$, the question is plausible. For the first 20 terms, the thing is possible I believe.

### Real calculation of all the digits of $\pi$

The process would go further by calculating all the figures. Indeed, by pushing the formula with many more terms, we quite easily reach a few million digits, but the speed of convergence is not geometric, it is polynomial even if the degree is very high. If we take the largest known value of B_2n (100 million), we will therefore have the series (or product if we take the Euler product) as follows.

$$\left( \frac{2n!}{B_n 2^n \left(1 - \frac{1}{2^n}\right)\left(1 - \frac{1}{3^n}\right)\left(1 - \frac{1}{5^n}\right)\left(1 - \frac{1}{7^n}\right)\cdots} \right)^{1/n}$$

So, with n = 100 million, and even if we take 1000 billion terms, we will have at best 1.2 billion digits of precision (1000 billion = $10^{12}$) so 12 digits only. This is a lot but little compared to the 62,000 billion digits obtained recently with the Chudnovsky formula, the latter has a convergence speed which is geometric, of the order of 14 decimal digits per term.

### Other asymptotic formulas containing the number π

First the factorial due to Stirling,

$$n! \approx \left(\frac{n}{e}\right)^n \sqrt{2\pi n} \left[1 + \frac{1}{12n} + \frac{1}{288\,n^2} - \frac{139}{51840\,n^3} + \cdots \right]$$

Which, once we isolate $\pi$ and use the bootstrap method,

$$\pi \approx \frac{\Gamma(n)^2\, n^{1-2n} e^{2n}}{2\left(1 + \frac{1}{6n}\right)\left(1 - \frac{5}{36n^2}\right)\left(1 + \frac{1}{72n^2}\right)}$$

Will give some digits of $\pi$ but even with n = $10^9$ only 30 digits of precision is obtained. It is possible to stretch the first expression as desired but this convergence is polynomial, the gain is marginal.

Another one used by Ramanujan and Hardy for the partition of n, $p(n)$, is approximately equal to

$$p(n) \approx \frac{1}{4n\sqrt{3}} e^{\pi\sqrt{2n/3}}$$

When π is isolated then we have

$$\pi \approx \frac{\ln(48\, p(n)^2 n^2)\sqrt{6}}{4\sqrt{n}}$$

Which also has a fairly slow convergence even if $n = 10^9$, here we could extend the expression with Rademacher's formula, which is much more complex but which makes the number π appear in several places, which makes the inversion operation for isolate rather difficult. This possibility cannot be ruled out at this time.

With the central binomial coefficient, the following approximation holds.

$$\binom{2n}{n} \approx \frac{4^n}{\sqrt{\pi n}}$$

Again, by isolating π we get

$$\pi \approx \frac{16^n}{n \binom{2n}{n}^2}$$

That can be pushed further by using the classic asymptotic expansion,

$$\pi \approx \frac{16^n}{n \binom{2n}{n}^2 \left[1 + \frac{1}{4n} + \frac{1}{32n^2} - \frac{1}{128n^3} - \frac{5}{2048n^4} \cdots\right]}$$

This time the approximation is the best so far, for $n = 10^9$ we get 100 digits of precision. The precision obtained if n ≫1 is much higher than the other 2 formulas but is not sufficient to obtain a precision of the order of the nth bit even if n is very large. If $n = 10^9$ we still obtain a hundred decimal places of precision. The asymptotic expansion is known and brings out the Bernoulli numbers which grow in size quite quickly beyond 14 terms. So the formulas obtained with Euler or Bernoulli numbers are by far the most accurate.